\documentclass[12pt]{article}%
\usepackage{amsmath}
\usepackage{amsfonts}
\usepackage{amssymb}
\usepackage{graphicx}
\usepackage{color}%
\providecommand{\U}[1]{\protect\rule{.1in}{.1in}}
\newtheorem{theorem}{Theorem}[section]

\newtheorem{lemma}[theorem]{Lemma}

\newtheorem{remark}[theorem]{Remark}

\newenvironment{proof}[1][Proof]{\textbf{#1.} }{\ \rule{0.5em}{0.5em}}

\makeatletter\@addtoreset{equation}{section}\makeatother
\evensidemargin=\oddsidemargin
\newdimen\dummy
\dummy=\oddsidemargin
\begin{document}

\title{An explicit factorization of the Green's function for an acoustic half-space
problem with impedance boundary conditions into an oscillatory exponential and
a slowly varying function.}
\author{C. Lin\thanks{(chuhe.lin@math.uzh.ch), Institut f\"{u}r Mathematik,
Universit\"{a}t Z\"{u}rich, Winterthurerstr 190, CH-8057 Z\"{u}rich,
Switzerland}
\and J.M. Melenk\thanks{(melenk@tuwien.ac.at), Institut f\"{u}r Analysis und
Scientific Computing, Technische Universit\"{a}t Wien, Wiedner Hauptstrasse
8-10, A-1040 Wien, Austria.}
\and S. Sauter\thanks{(stas@math.uzh.ch), Institut f\"{u}r Mathematik,
Universit\"{a}t Z\"{u}rich, Winterthurerstr 190, CH-8057 Z\"{u}rich,
Switzerland}}
\maketitle

\begin{abstract}
In this paper, new representations of the Green's function for an acoustic
$d$-dimensional half-space problem with impedance boundary conditions are
presented. The main features of the new representation are:

a) in addition to additive terms that appear also in the case of Dirichlet or
Neumann boundary conditions, the remaining part of the Green's function is
factored into an oscillatory complex exponential function (with the product of
the wavenumber and the eikonal as argument) and a remaining function which is
slowly varying and hence allows for efficient polynomial approximation;

b) the representation is given uniformly for all parameters by a single
formula which consists of the product of two analytic functions.

\end{abstract}



\section{Introduction}

In this paper we consider an acoustic half-space problem with impedance
boundary conditions in general $d$ spatial dimensions. The main result is the
derivation of a new integral representation of the corresponding Green's
function in a form where oscillatory Fourier-type integrals (see, e.g.,
\cite[(13)]{ChandlerWilde_imp_Green}, \cite[(21)]{Hein_Nedelec_green},
\cite{Duran_stat_phase}, \cite{Duran_stat_phase_long}, \cite{hoernig2010green}%
, \cite{ochmann2011closed}, \cite{Gimperlein_halfspace}, \cite{ROJAS20183903})
are avoided so that it is well suited for an analysis, the study of its
approximation, and the derivation of uniform (high order) asymptotic
expansions. In contrast to the representations cited above, the integrand in
the new integral representation is non-oscillatory with respect to the outer
variable and defines a function which is non-oscillatory. For odd spatial
degree and impedance parameter $\beta=1$, i.e., the half-space problem with
Robin boundary conditions we present fully explicit representations of this
Green's function.

While the focus in this paper is on the derivation of the new representations,
the companion paper \cite{LinMelenkSauter_Gimp_II} is devoted to its efficient approximation.

\section{The acoustic half-space problem with impedance boundary conditions}

Let the upper half-space in $\mathbb{R}^{d}$, $d\in\left\{  1,2,\ldots
\right\}  $, and its boundary be denoted by%
\begin{align*}
H_{+} &  :=\left\{  \mathbf{x}=\left(  x_{j}\right)  _{j=1}^{d}\in
\mathbb{R}^{d}\mid x_{d}>0\right\}  ,\\
H_{0} &  :=\partial H_{+}:=\left\{  \mathbf{x}=\left(  x_{j}\right)
_{j=1}^{d}\in\mathbb{R}^{d}\mid x_{d}=0\right\}  .
\end{align*}
The outward normal vector is given $\mathbf{n}=\left(  0,\ldots,0,-1\right)
^{T}$. Let%
\[
\mathbb{C}_{>0}:=\left\{  \zeta\in\mathbb{C}\mid\operatorname{Re}%
\zeta>0\right\}  .
\]

We consider the problem to find the Green's function $G:H_{+}\times
H_{+}\rightarrow\mathbb{C}$ for the acoustic half-plane problem with impedance
boundary conditions:%
\begin{equation}%
\begin{array}
[c]{rll}%
-\Delta_{\mathbf{x}}G\left(  \mathbf{x},\mathbf{y}\right)  +s^{2}G\left(
\mathbf{x},\mathbf{y}\right)   & =\delta_{0}\left(  \mathbf{x}-\mathbf{y}%
\right)   & \text{for }\left(  \mathbf{x},\mathbf{y}\right)  \in H_{+}\times
H_{+},\\
\frac{\partial}{\partial\mathbf{n}_{\mathbf{x}}}G\left(  \mathbf{x}%
,\mathbf{y}\right)  +s\beta G\left(  \mathbf{x},\mathbf{y}\right)   & =0 &
\text{for }\left(  \mathbf{x},\mathbf{y}\right)  \in H_{0}\times H_{+},\\
G\left(  r%
\mbox{\boldmath$ \zeta$}%
,\mathbf{y}\right)   & \overset{r\rightarrow+\infty}{\rightarrow}0 & \text{for
}\left(
\mbox{\boldmath$ \zeta$}%
,\mathbf{y}\right)  \in H_{+}\times H_{+}%
\end{array}
\label{goveq}%
\end{equation}
for some $\beta>0$ and frequency $s\in\mathbb{C}_{>0}$. The index $\mathbf{x}$
in the differential operators indicate that the derivative is taken with
respect to the variable $\mathbf{x}$.

\begin{remark}
\label{RemLimitAbsorption}Problem (\ref{goveq}) is formulated for
$s\in\mathbb{C}_{>0}$. The Green's function $G=G_{s}$ depends on $s$ and for
$\operatorname{Re}s>0$ it is assumed to decay for $\mathbf{x}=r%
\mbox{\boldmath$ \zeta$}%
$ as $r\rightarrow+\infty$ for any fixed direction $%
\mbox{\boldmath$ \zeta$}%
\in H_{+}$. Problem (\ref{goveq}) for the case $s\in\operatorname*{i}%
\mathbb{R}\backslash\left\{  0\right\}  $ is considered as the limit case from
the positive complex half-plane $\mathbb{C}_{>0}$:%
\[
G_{s}=\lim_{\substack{\sigma\rightarrow s\\\sigma\in\mathbb{C}_{>0}}%
}G_{\sigma}\text{.}%
\]

\end{remark}

The case $d=1$ can be solved fully explicitly.

\begin{remark}
For $d=1$, the Green's function for the impedance problem (\ref{goveq}) is
given by%
\begin{equation}
G\left(  x,y\right)  =\frac{1}{2s}\left(  \operatorname*{e}%
\nolimits^{-s\left\vert x-y\right\vert }+\frac{1-\beta}{1+\beta}%
\operatorname*{e}\nolimits^{-s\left(  x+y\right)  }\right)  .\label{qdodd1}%
\end{equation}

\end{remark}

For the rest of the paper we assume that $d\in\left\{  2,3,\ldots\right\}  $
and introduce%
\[
\nu:=\left(  d-3\right)  /2.
\]

\section{Representation of the Green's function}

The representation of the Green's function as the solution of (\ref{goveq})
requires some preparations. Let $K_{\nu}$ denote the Macdonald function
(modified Bessel function of the second kind and order $\nu$, see, e.g.,
\cite[\S 10.25]{NIST:DLMF}, \cite{Macdonald_fct_orginal}). We introduce the
function%
\begin{equation}
g_{\nu}\left(  r\right)  :=\frac{1}{\left(  2\pi\right)  ^{\nu+3/2}}\left(
\frac{s}{r}\right)  ^{\nu+1/2}K_{\nu+1/2}\left(  sr\right)  \label{defgnue}%
\end{equation}
and note that $g_{\nu}\left(  \left\Vert \mathbf{x}-\mathbf{y}\right\Vert
\right)  $ is the full space Green's function for the Helmholtz operator (see
\cite[(9.14)]{Mclean00} and \cite[(6), (12)]{Buslaev_1975} in combination with
the connecting formula \cite[\S 10.27.8]{NIST:DLMF}). For $\mathbf{y=}\left(
y_{j}\right)  _{j=1}^{d}\in H_{+}$, we introduce the reflection operator
$\mathbf{Ry}=\left(  \mathbf{y}^{\prime},-y_{d}\right)  $, where
$\mathbf{y}^{\prime}=\left(  y_{j}\right)  _{j=1}^{d-1}$.

Let the functions $r:\mathbb{R}^{d}\rightarrow\mathbb{R}$ and $r_{+}%
:\mathbb{R}^{d}\rightarrow\mathbb{R}$ be defined for $\mathbf{z}\in H_{+}$ and
$\mathbf{z}^{\prime}:=\left(  z_{j}\right)  _{j=1}^{d-1}$ by%

\[
r\left(  \mathbf{z}\right)  :=\left\Vert \mathbf{z}\right\Vert ,\quad
r_{+}\left(  \mathbf{z}\right)  :=r\left(  \mathbf{z}\right)  +\beta z_{d}%
\]
and set%
\begin{equation}
y\left(  \mathbf{z,\cdot}\right)  :\left[  z_{d},\infty\right[  \rightarrow
\left[  0,\infty\right[  ,\quad y\left(  \mathbf{z},t\right)  :=-r_{+}\left(
\mathbf{z}\right)  +\beta t+\mu\left(  \mathbf{z}^{\prime},t\right)
\label{defyvont}%
\end{equation}
with the function $\mu\left(  \mathbf{z}^{\prime},\cdot\right)  :\left[
z_{d},\infty\right[  \rightarrow\left[  \left\Vert \mathbf{z}\right\Vert
,\infty\right[  $ given by%
\[
\mu\left(  \mathbf{z}^{\prime},t\right)  :=\sqrt{\left\Vert \mathbf{z}%
^{\prime}\right\Vert ^{2}+t^{2}}.
\]
The derivative of $y$ satisfies%
\begin{equation}
\frac{\partial y\left(  \mathbf{z},t\right)  }{\partial t}=\beta+\frac{t}%
{\mu\left(  \mathbf{z}^{\prime},t\right)  }>0 \label{cfdydt}%
\end{equation}
so that $y\left(  \mathbf{z},\cdot\right)  $ maps the interval $\left[
z_{d},\infty\right[  $ strictly increasing onto $\left[  0,\infty\right[  $.
Its inverse%
\begin{equation}
t\left(  \mathbf{z},\cdot\right)  :\left[  0,\infty\right[  \rightarrow\left[
z_{d},\infty\right[  \label{tsub}%
\end{equation}
is also strictly increasing. The derivative $\partial t\left(  \mathbf{z}%
,y\right)  /\partial y$ can be expressed by using (\ref{cfdydt}):%
\begin{equation}
\frac{\partial t\left(  \mathbf{z},y\right)  }{\partial y}=\frac{\tilde{\mu
}\left(  \mathbf{z},y\right)  }{t\left(  \mathbf{z},y\right)  +\beta\tilde
{\mu}\left(  \mathbf{z},y\right)  }, \label{dert}%
\end{equation}
where%
\begin{equation}
\tilde{\mu}\left(  \mathbf{z},y\right)  :=\mu\left(  \mathbf{z}^{\prime
},t\left(  \mathbf{z},y\right)  \right)  \text{\quad and\quad}\frac
{\partial\tilde{\mu}\left(  \mathbf{z},y\right)  }{\partial y}=\frac{t\left(
\mathbf{z},y\right)  }{t\left(  \mathbf{z},y\right)  +\beta\tilde{\mu}\left(
\mathbf{z},y\right)  }>0. \label{defsmue}%
\end{equation}
In the following, the shorthands%
\begin{equation}
r=r\left(  \mathbf{z}\right)  ,\quad t=t\left(  \mathbf{z},y\right)
,\quad\tilde{\mu}=\tilde{\mu}\left(  \mathbf{z},y\right)  \label{shorthands}%
\end{equation}
will be used. A key role for the representation of the Green's function will
be played by the function%
\begin{equation}
\psi_{\nu,s}\left(  \mathbf{z}\right)  :=\int_{0}^{\infty}\frac
{\operatorname*{e}\nolimits^{-sy}}{t+\beta\tilde{\mu}}\frac{\operatorname*{e}%
\nolimits^{s\tilde{\mu}}K_{\nu+1/2}\left(  s\tilde{\mu}\right)  }{\left(
s\tilde{\mu}\right)  ^{\nu-1/2}}dy. \label{defQ1nue}%
\end{equation}
This integral exists as an improper Riemann integral for any $s\in
\mathbb{C}_{>0}$ while for $s\in\operatorname*{i}\mathbb{R}\backslash\left\{
0\right\}  $ it is defined as the limit described in Remark
\ref{RemLimitAbsorption} and discussed in more detail in Remark
\ref{Rempartint}.

\begin{theorem}
\label{TheoMain}Let $d\in\left\{  2,3,\ldots\right\}  $ denote the spatial
dimension. The Green's function for the acoustic half-space problem with
impedance boundary conditions is
\[
G\left(  \mathbf{x},\mathbf{y}\right)  =g_{\nu}\left(  \left\Vert
\mathbf{x}-\mathbf{y}\right\Vert \right)  +g_{\nu}\left(  \left\Vert
\mathbf{x}-\mathbf{Ry}\right\Vert \right)  +G_{\operatorname*{imp}}\left(
\mathbf{x}-\mathbf{Ry}\right)  ,
\]
where $g_{\nu}$ is as in (\ref{defgnue}), and%
\[
G_{\operatorname*{imp}}\left(  \mathbf{z}\right)  =-\frac{\beta}{\pi}\left(
\frac{s^{2}}{2\pi}\right)  ^{\nu+1/2}\operatorname*{e}\nolimits^{-s\left\Vert
\mathbf{z}\right\Vert }\psi_{\nu,s}\left(  \mathbf{z}\right)
\]
solves the governing equation (\ref{goveq}).
\end{theorem}

%

\proof
Let $s\in\mathbb{C}_{>0}$ and let the Helmholtz operator be denoted by
$\mathcal{L}_{s}=-\Delta+s^{2}$. Since $g_{\nu}$ is the full space Green's
function it holds $\mathcal{L}_{\mathbf{x},s}g_{\nu}\left(  \left\Vert
\mathbf{x}-\mathbf{y}\right\Vert \right)  =\delta_{0}\left(  \mathbf{x}%
-\mathbf{y}\right)  $ for all $\mathbf{x},\mathbf{y}\in H_{+}$ in a distributional sense. It follows via
the chain rule $\mathcal{L}_{\mathbf{x},s}g_{\nu}\left(  \left\Vert
\mathbf{x}-\mathbf{Ry}\right\Vert \right)  =0$ for all $\mathbf{x}%
,\mathbf{y}\in H_{+}$. Let $\mathbf{z}:=\mathbf{x}-\mathbf{Ry}$,
$r:=\left\Vert \mathbf{z}\right\Vert $, and note that $\mathbf{z}\in H_{+}$.
To show that $G_{\operatorname*{imp}}$ is Helmholtz-harmonic we apply the
variable transform (cf. (\ref{defyvont})) $y\leftarrow y\left(
\mathbf{z,\cdot}\right)  $ and obtain%
\begin{equation}
G_{\operatorname*{imp}}\left(  \mathbf{z}\right)  =-2\beta\left(  \frac
{s}{2\pi}\right)  ^{\nu+3/2}\int_{z_{d}}^{\infty}\operatorname*{e}%
\nolimits^{-s\beta\left(  t-z_{d}\right)  }\frac{K_{\nu+1/2}\left(
s\mu\left(  \mathbf{z}^{\prime},t\right)  \right)  }{\left(  \mu\left(
\mathbf{z}^{\prime},t\right)  \right)  ^{\nu+1/2}}dt. \label{Gimphandy}%
\end{equation}
The derivatives with respect to $z_{d}$ are given by (using $\left(
z^{-\lambda}K_{\lambda}\left(  z\right)  \right)  ^{\prime}=-z^{-\lambda
}K_{\lambda+1}\left(  z\right)  $; see \cite[10.29.4]{NIST:DLMF})%
\begin{align}
\frac{\partial G_{\operatorname*{imp}}\left(  \mathbf{z}\right)  }{\partial
z_{d}}  &  =2\beta\left(  \frac{s}{2\pi}\right)  ^{\nu+3/2}\frac{K_{\nu
+1/2}\left(  sr\right)  }{r^{\nu+1/2}}+s\beta G_{\operatorname*{imp}}\left(
\mathbf{z}\right) \label{firstzdderivative}\\
\frac{\partial^{2}G_{\operatorname*{imp}}\left(  \mathbf{z}\right)  }{\partial
z_{d}^{2}}  &  =-2s\beta z_{d}\left(  \frac{s}{2\pi}\right)  ^{\nu+3/2}%
\frac{K_{\nu+3/2}\left(  sr\right)  }{r^{\nu+3/2}}+2s\beta^{2}\left(  \frac
{s}{2\pi}\right)  ^{\nu+3/2}\frac{K_{\nu+1/2}\left(  sr\right)  }{r^{\nu+1/2}%
}+s^{2}\beta^{2}G_{\operatorname*{imp}}\left(  \mathbf{z}\right)  .\nonumber
\end{align}
To simplify this expression we apply two times integration by parts to
$G_{\operatorname*{imp}}$ as in (\ref{Gimphandy})%
\begin{align*}
G_{\operatorname*{imp}}  &  =-\frac{2}{s}\left(  \frac{s}{2\pi}\right)
^{\nu+3/2}\frac{K_{\nu+1/2}\left(  sr\right)  }{r^{\nu+1/2}}+2\left(  \frac
{s}{2\pi}\right)  ^{\nu+3/2}\int_{z_{d}}^{\infty}\operatorname*{e}%
\nolimits^{-s\beta\left(  t-z_{d}\right)  }\frac{tK_{\nu+3/2}\left(
s\mu\left(  \mathbf{z}^{\prime},t\right)  \right)  }{\left(  \mu\left(
\mathbf{z}^{\prime},t\right)  \right)  ^{\nu+3/2}}dt\\
&  =-\frac{2}{s}\left(  \frac{s}{2\pi}\right)  ^{\nu+3/2}\frac{K_{\nu
+1/2}\left(  sr\right)  }{r^{\nu+1/2}}+\frac{2z_{d}}{s\beta}\left(  \frac
{s}{2\pi}\right)  ^{\nu+3/2}\frac{K_{\nu+3/2}\left(  sr\right)  }{r^{\nu+3/2}%
}\\
&  +\frac{2}{s\beta}\left(  \frac{s}{2\pi}\right)  ^{\nu+3/2}\int_{z_{d}%
}^{\infty}\operatorname*{e}\nolimits^{-s\beta\left(  t-z_{d}\right)  }\left(
\frac{K_{\nu+3/2}\left(  s\mu\left(  \mathbf{z}^{\prime},t\right)  \right)
}{\left(  \mu\left(  \mathbf{z}^{\prime},t\right)  \right)  ^{\nu+3/2}}%
-\frac{st^{2}K_{\nu+5/2}\left(  s\mu\left(  \mathbf{z}^{\prime},t\right)
\right)  }{\left(  \mu\left(  \mathbf{z}^{\prime},t\right)  \right)
^{\nu+5/2}}\right)  dt.
\end{align*}
In this way, we get with $t^{2}=\mu^{2}-\left\Vert \mathbf{z}^{\prime
}\right\Vert ^{2}$%
\begin{align}
\frac{\partial^{2}G_{\operatorname*{imp}}\left(  \mathbf{z}\right)  }{\partial
z_{d}^{2}}  &  =2s\beta\left(  \frac{s}{2\pi}\right)  ^{\nu+3/2}\int_{z_{d}%
}^{\infty}\operatorname*{e}\nolimits^{-s\beta\left(  t-z_{d}\right)  }%
\times\nonumber\\
&  \times\left(  \frac{K_{\nu+3/2}\left(  s\mu\left(  \mathbf{z}^{\prime
},t\right)  \right)  }{\left(  \mu\left(  \mathbf{z}^{\prime},t\right)
\right)  ^{\nu+3/2}}-s\frac{K_{\nu+5/2}\left(  s\mu\left(  \mathbf{z}^{\prime
},t\right)  \right)  }{\left(  \mu\left(  \mathbf{z}^{\prime},t\right)
\right)  ^{\nu+1/2}}+s\left\Vert \mathbf{z}^{\prime}\right\Vert ^{2}%
\frac{K_{\nu+5/2}\left(  s\mu\left(  \mathbf{z}^{\prime},t\right)  \right)
}{\left(  \mu\left(  \mathbf{z}^{\prime},t\right)  \right)  ^{\nu+5/2}%
}\right)  dt. \label{Gimpderzd2}%
\end{align}
For the gradient and the Laplacian with respect to $\mathbf{z}^{\prime}$, we
calculate%
\begin{align*}
\nabla_{\mathbf{z}^{\prime}}\frac{K_{\nu+1/2}\left(  s\mu\left(
\mathbf{z}^{\prime},t\right)  \right)  }{\left(  \mu\left(  \mathbf{z}%
^{\prime},t\right)  \right)  ^{\nu+1/2}}  &  =-s\mathbf{z}^{\prime}%
\frac{K_{\nu+3/2}\left(  s\mu\left(  \mathbf{z}^{\prime},t\right)  \right)
}{\left(  \mu\left(  \mathbf{z}^{\prime},t\right)  \right)  ^{\nu+3/2}},\\
\Delta_{\mathbf{z}^{\prime}}\frac{K_{\nu+1/2}\left(  s\mu\left(
\mathbf{z}^{\prime},t\right)  \right)  }{\left(  \mu\left(  \mathbf{z}%
^{\prime},t\right)  \right)  ^{\nu+1/2}}  &  =-2s\left(  \nu+1\right)
\frac{K_{\nu+3/2}\left(  s\mu\left(  \mathbf{z}^{\prime},t\right)  \right)
}{\left(  \mu\left(  \mathbf{z}^{\prime},t\right)  \right)  ^{\nu+3/2}}%
+s^{2}\left\Vert \mathbf{z}^{\prime}\right\Vert ^{2}\frac{K_{\nu+5/2}\left(
s\mu\left(  \mathbf{z}^{\prime},t\right)  \right)  }{\left(  \mu\left(
\mathbf{z}^{\prime},t\right)  \right)  ^{\nu+5/2}}.
\end{align*}
Then we combine this with (\ref{Gimpderzd2}) to obtain%
\begin{align*}
-\Delta G_{\operatorname*{imp}}\left(  \mathbf{z}\right)  =-2s^{2}\beta\left(
\frac{s}{2\pi}\right)  ^{\nu+3/2}\int_{z_{d}}^{\infty}  &  \frac
{\operatorname*{e}\nolimits^{-s\beta\left(  t-z_{d}\right)  }}{\left(
\mu\left(  \mathbf{z}^{\prime},t\right)  \right)  ^{\nu+1/2}}\times\\
&  \times\left(  \left(  2\nu+3\right)  \frac{K_{\nu+3/2}\left(  s\mu\left(
\mathbf{z}^{\prime},t\right)  \right)  }{s\mu\left(  \mathbf{z}^{\prime
},t\right)  }-K_{\nu+5/2}\left(  s\mu\left(  \mathbf{z}^{\prime},t\right)
\right)  \right)  dt.
\end{align*}
Next we use \cite[10.29.1]{NIST:DLMF}, i.e., $K_{\lambda+1}\left(  z\right)
=\frac{2\lambda}{z}K_{\lambda}\left(  z\right)  +K_{\lambda-1}\left(
z\right)  $ for $\lambda=\nu+3/2$, and get%
\[
-\Delta G_{\operatorname*{imp}}\left(  \mathbf{z}\right)  =s^{2}2\beta\left(
\frac{s}{2\pi}\right)  ^{\nu+3/2}\int_{z_{d}}^{\infty}\operatorname*{e}%
\nolimits^{-s\beta\left(  t-z_{d}\right)  }\frac{K_{\nu+1/2}\left(
s\mu\left(  \mathbf{z}^{\prime},t\right)  \right)  }{\left(  \mu\left(
\mathbf{z}^{\prime},t\right)  \right)  ^{\nu+1/2}}dt=-s^{2}%
G_{\operatorname*{imp}}\left(  \mathbf{z}\right)  .
\]
This implies $\mathcal{L}_{\mathbf{x},s}G_{\operatorname*{imp}}\left(
\mathbf{x-Ry}\right)  =0$ and, in turn, $\mathcal{L}_{\mathbf{x},s}G\left(
\mathbf{x},\mathbf{y}\right)  =\delta_{0}\left(  \mathbf{x}-\mathbf{y}\right)
$.

To verify the boundary condition we denote by $\mathcal{B}_{\mathbf{x}%
,s}:=\partial/\partial\mathbf{n}_{\mathbf{x}}+s\beta$ the boundary
differential operator in (\ref{goveq}). Since $\left.  \left\Vert
\mathbf{x-y}\right\Vert \right\vert _{x_{d}=0}=\left.  \left\Vert \left(
\mathbf{x-Ry}\right)  \right\Vert \right\vert _{x_{d}=0}$ we obtain%
\begin{equation}
\mathcal{B}_{\mathbf{x},s}G\left(  \mathbf{x},\mathbf{y}\right)  =\left.
2s\beta g_{\nu}\left(  \left\Vert \mathbf{x}-\mathbf{y}\right\Vert \right)
\right\vert _{x_{d}=0}+\mathcal{B}_{\mathbf{x},s}G_{\operatorname*{imp}%
}\left(  \mathbf{x}-\mathbf{Ry}\right)  .\label{boundop1}%
\end{equation}
From (\ref{firstzdderivative}) it follows that the normal derivative of
$G_{\operatorname*{imp}}$ has the form%
\[
\frac{\partial}{\partial\mathbf{n}_{\mathbf{z}}}G_{\operatorname*{imp}}\left(
\mathbf{z}\right)  =-\frac{\partial}{\partial z_{d}}G_{\operatorname*{imp}%
}\left(  \mathbf{z}\right)  =-2\beta\left(  \frac{s}{2\pi}\right)  ^{\nu
+3/2}\frac{K_{\nu+1/2}\left(  sr\right)  }{r^{\nu+1/2}}-s\beta
G_{\operatorname*{imp}}\left(  \mathbf{z}\right)  .
\]
We use (\ref{defgnue}) for the second equality in%
\[
\mathcal{B}_{\mathbf{x},s}G_{\operatorname*{imp}}\left(  \mathbf{x}%
-\mathbf{Ry}\right)  =\left.  -2\beta\left(  \frac{s}{2\pi}\right)  ^{\nu
+3/2}\frac{K_{\nu+1/2}\left(  s\left\Vert \mathbf{x}-\mathbf{Ry}\right\Vert
\right)  }{\left\Vert \mathbf{x}-\mathbf{Ry}\right\Vert ^{\nu+1/2}}\right\vert
_{x_{d}=0}=-2s\beta\left.  g_{\nu}\left(  \left\Vert \mathbf{x}-\mathbf{y}%
\right\Vert \right)  \right\vert _{x_{d}=0},
\]
and a comparison with (\ref{boundop1}) leads to $\mathcal{B}_{\mathbf{x}%
,s}G\left(  \mathbf{x},\mathbf{y}\right)  =0$ {for $\mathbf{x}\in H_{0}$.}

Finally, we investigate the decay condition and recall $\operatorname{Re}s>0$.
The asymptotics for modified Bessel functions for large argument are
well-known (see, e.g., \cite[10.40.2]{NIST:DLMF}) to be $K_{\nu}\left(
z\right)  \sim\sqrt{\frac{\pi}{2z}}\operatorname*{e}\nolimits^{-z}$. This
directly implies the decay of $g_{\nu}\left(  \left\Vert \mathbf{x}%
-\mathbf{y}\right\Vert \right)  +g_{\nu}\left(  \left\Vert \mathbf{x}%
-\mathbf{Ry}\right\Vert \right)  $.

For $G_{\operatorname*{imp}}$ we start from (\ref{Gimphandy}) and estimate%
\[
\left\vert G_{\operatorname*{imp}}\left(  \mathbf{z}\right)  \right\vert
\leq2\beta\left(  \frac{\left\vert s\right\vert }{2\pi}\right)  ^{\nu
+3/2}M_{\nu}\left(  \mathbf{z}\right)  \int_{z_{d}}^{\infty}\operatorname*{e}%
\nolimits^{-(  \operatorname{Re}(s))  \ \beta\left(  t-z_{d}\right)
\ }dt=2\beta\left(  \frac{\left\vert s\right\vert }{2\pi}\right)  ^{\nu
+3/2}\frac{M_{\nu}\left(  \mathbf{z}\right)  }{\beta\operatorname{Re}(s)}.
\]
with%
\[
M_{\nu}\left(  \mathbf{z}\right)  :=\sup_{t\in\left]  z_{d},\infty\right[
}\left\vert \frac{K_{\nu+1/2}\left(  s\mu\left(  \mathbf{z}^{\prime},t\right)
\right)  }{\left(  \mu\left(  \mathbf{z}^{\prime},t\right)  \right)
^{\nu+1/2}}\right\vert .
\]
Note that $\left\vert \mu\left(  \mathbf{z}^{\prime},t\right)  \right\vert
\geq\left\Vert \mathbf{z}\right\Vert $ so that the exponential decay of
$K_{\nu+1/2}$ for large argument implies the decay of $G_{\operatorname*{imp}%
}$ and, in turn, of $G$ as required in the third condition in (\ref{goveq}).%
\endproof

\begin{remark}
\label{Rempartint}For $d=2,3$ and $s\in\operatorname*{i}\mathbb{R}%
\backslash\left\{  0\right\}  $, the integral (\ref{defQ1nue}) does not exist
as an improper Riemann integral, and the limit $\lim_{\substack{\sigma
\rightarrow s\\\sigma\in\mathbb{C}_{>0}}}G_{\sigma}$ cannot be interchanged
with the integral. This problem can be resolved by an integration by parts: we
set (with shorthands (\ref{shorthands}))%
\[
q_{\nu}\left(  \mathbf{z},y\right)  :=\frac{d}{dy}\left(  \frac
{\operatorname*{e}\nolimits^{s\tilde{\mu}}K_{\nu+1/2}\left(  s\tilde{\mu
}\right)  }{\left(  t+\beta\tilde{\mu}\right)  \left(  s\tilde{\mu}\right)
^{\nu-1/2}}\right)
\]
and obtain%
\begin{equation}
\psi_{\nu,s}\left(  \mathbf{z}\right)  =\frac{\operatorname*{e}^{sr}%
K_{\nu+1/2}\left(  sr\right)  }{s\left(  z_{d}+\beta r\right)  \left(
sr\right)  ^{\nu-1/2}}+\frac{1}{s}\int_{0}^{\infty}\operatorname*{e}%
\nolimits^{-sy}q_{\nu}\left(  \mathbf{z},y\right)  dy. \label{Defpsinues}%
\end{equation}
The integral in (\ref{Defpsinues}) exists as an improper Riemann integral also
for the cases $d=2,3$ and $s\in\operatorname*{i}\mathbb{R}\backslash\left\{
0\right\}  $ (see \cite[Lem. 4.11(1)]{LinMelenkSauter_Gimp_II}).
\end{remark}

\section{Special cases}

If the space dimension $d$ is odd, i.e., $\nu$ is an integer, and $\beta=1$,
the functions $G_{\operatorname*{imp}}\left(  \mathbf{z}\right)  $ allow for a
more explicit representation compared to the integral in (\ref{defQ1nue}).

\begin{lemma}
Let $d\in\left\{  3,5,\ldots\right\}  $ so that the parameter $\nu=\left(
d-3\right)  /2$ is an integer and assume $\beta=1$. In this case, the function
$G_{\operatorname*{imp}}\left(  \mathbf{z}\right)  $ is given for $d=3$, i.e.,
$\nu=0$, by%
\begin{subequations}
\label{qdodd}
\end{subequations}%
\begin{equation}
G_{\operatorname*{imp}}\left(  \mathbf{z}\right)  =-\left(  \frac{s}{2\pi
}\right)  \operatorname*{e}\nolimits^{-s\left\Vert \mathbf{z}\right\Vert
}U\left(  1,1,s\left(  \left\Vert \mathbf{z}\right\Vert +z_{3}\right)
\right)  \tag{%
\ref{qdodd}%
a}\label{qdodda}%
\end{equation}
with Tricomi's (confluent hypergeometric) function $U\left(  a,b,z\right)  $
(other name: Gordon function), (see \cite[13.2.6]{NIST:DLMF},
\cite{Tricomi_Kummer}, \cite[p. 671]{MorseFeshbach}), which is a solution of
Kummer's differential equation (see \cite[(9.)]{Kummer_original}).

For $\nu=1,2,3,\ldots$ it holds%
\begin{equation}
G_{\operatorname*{imp}}\left(  \mathbf{z}\right)  =\left(  s-\frac{\partial
}{\partial z_{d}}\right)  ^{\nu-1}\Psi_{\nu,s}\left(  \mathbf{z}\right)  \tag{%
\ref{qdodd}%
b}\label{qdoddc}%
\end{equation}
with $\Psi_{\nu,s}$ defined by%
\begin{equation}
\Psi_{\nu,s}\left(  \mathbf{z}\right)  :=-\frac{s}{\left(  2\pi\right)
^{\nu+1}}\frac{\operatorname*{e}\nolimits^{-s\left\Vert \mathbf{z}\right\Vert
}}{\left(  \left\Vert \mathbf{z}\right\Vert +z_{d}\right)  ^{\nu}\left\Vert
\mathbf{z}\right\Vert }. \label{Psinues}%
\end{equation}

\end{lemma}

%

\proof
We prove this lemma for $s\in\mathbb{C}_{>0}$, while the case $s\in
\operatorname*{i}\mathbb{R}\backslash\left\{  0\right\}  $ is obtained by
taking the limit $\tilde{s}\rightarrow s$ from $\mathbb{C}_{>0}$ in
(\ref{qdodd}) and (\ref{Psinues}).

For $\mathbf{x},\mathbf{y}\in H_{+}$, the notation\quad%
\begin{equation}
\mathbf{z}:=\mathbf{x}-\mathbf{Ry}\text{,\quad}r=\left\Vert \mathbf{z}%
\right\Vert \text{,\quad}\mathbf{z}^{\prime}=\left(  z_{j}\right)
_{j=1}^{d-1}\text{,\quad}\omega:=\left\Vert \mathbf{z}^{\prime}\right\Vert
,\quad\text{and }\tilde{q}\left(  x\right)  :=\sqrt{x^{2}+s^{2}}
\label{shorthand2}%
\end{equation}
is used. Note that $K_{1/2}\left(  z\right)  =\sqrt{\pi/\left(  2z\right)
}\operatorname*{e}^{-z}$ (cf., e.g., \cite[(5)]{Bessel_polys}). The
representation (\ref{Gimphandy}) of $G_{\operatorname*{imp}}$ for $\beta=1$
and $d=3$, i.e., $\nu=0$ takes the form%
\[
G_{\operatorname*{imp}}\left(  \mathbf{z}\right)  =-2\left(  \frac{s}{2\pi
}\right)  ^{3/2}\int_{z_{3}}^{\infty}\operatorname*{e}\nolimits^{-s\left(
t-z_{3}\right)  }\frac{K_{1/2}\left(  s\mu\left(  \mathbf{z}^{\prime
},t\right)  \right)  }{\left(  \mu\left(  \mathbf{z}^{\prime},t\right)
\right)  ^{1/2}}dt=-\frac{s}{2\pi}\int_{z_{3}}^{\infty}\frac{\operatorname*{e}%
\nolimits^{-s\left(  t-z_{3}+\mu\left(  \mathbf{z}^{\prime},t\right)  \right)
}}{\mu\left(  \mathbf{z}^{\prime},t\right)  }dt.
\]
The change of variables%
\[
y=t-z_{3}+\mu\left(  \mathbf{z}^{\prime},t\right)  -\left\Vert \mathbf{z}%
\right\Vert \quad\text{with\quad}\frac{dt}{dy}=1/\left(  dy/dt\right)
=\frac{\mu\left(  \mathbf{z}^{\prime},t\right)  }{\mu\left(  \mathbf{z}%
^{\prime},t\right)  +t}=\frac{\mu\left(  \mathbf{z}^{\prime},t\right)
}{\left\Vert \mathbf{z}\right\Vert +z_{3}+y}%
\]
leads to (with the exponential integral $\operatorname*{Ei}$; see
\cite[6.2.5]{NIST:DLMF})%
\begin{align*}
G_{\operatorname*{imp}}\left(  \mathbf{z}\right)   &  =-\frac
{s\operatorname*{e}\nolimits^{-s\left\Vert \mathbf{z}\right\Vert }}{2\pi}%
\int_{0}^{\infty}\frac{\operatorname*{e}\nolimits^{-sy}}{\left\Vert
\mathbf{z}\right\Vert +z_{3}+y}dy\overset{\text{\cite[3.352(2)]{gradstein}%
}}{=}\frac{s\operatorname*{e}\nolimits^{sz_{3}}}{2\pi}\operatorname{Ei}\left(
-s\left(  \left\Vert \mathbf{z}\right\Vert +z_{3}\right)  \right) \\
&  \overset{\text{\cite[6.2.6 \& 13.6.6]{NIST:DLMF}}}{=}-\frac
{s\operatorname*{e}\nolimits^{-s\left\Vert \mathbf{z}\right\Vert }}{2\pi
}U\left(  1,1,s\left(  \left\Vert \mathbf{z}\right\Vert +z_{3}\right)
\right)  .
\end{align*}

For the second claim (\ref{qdoddc}), it suffices to prove that
(\ref{Gimphandy}) for $\beta=1$ defines the same function as defined in
(\ref{qdoddc}). The function in (\ref{qdoddc}) is denoted by $\tilde
{G}_{\operatorname*{imp}}$ and the one in (\ref{Gimphandy}) for $\beta=1$ by
$G_{\operatorname*{imp}}$ so that the claim is $G_{\operatorname*{imp}}%
=\tilde{G}_{\operatorname*{imp}}$. The relation in (\ref{firstzdderivative})
implies that $G_{\operatorname*{imp}}$ satisfies the differential equation%
\begin{equation}
\left(  \frac{\partial}{\partial z_{d}}-s\right)  G_{\operatorname*{imp}%
}\left(  \mathbf{z}\right)  =2\left(  \frac{s}{2\pi}\right)  ^{\nu+3/2}%
\frac{K_{\nu+1/2}\left(  sr\right)  }{r^{\nu+1/2}} \label{ODE1}%
\end{equation}
and $G_{\operatorname*{imp}}\left(  \mathbf{z}\right)  $ decays to zero for
$z_{d}\rightarrow\infty$ as shown in the last part of the proof of Theorem
\ref{TheoMain}. Hence, it is sufficient to prove that $\tilde{G}%
_{\operatorname*{imp}}$ satisfies (\ref{ODE1}) and the decay condition.
Plugging in (\ref{qdoddc}) into (\ref{ODE1}) leads to the condition%
\begin{equation}
\left(  s-\frac{\partial}{\partial z_{d}}\right)  ^{\nu}\frac{s}{\left(
2\pi\right)  ^{\nu+1}}\frac{\operatorname*{e}\nolimits^{-s\left\Vert
\mathbf{z}\right\Vert }}{\left(  \left\Vert \mathbf{z}\right\Vert
+z_{d}\right)  ^{\nu}\left\Vert \mathbf{z}\right\Vert }=2\left(  \frac{s}%
{2\pi}\right)  ^{\nu+3/2}\frac{K_{\nu+1/2}\left(  sr\right)  }{r^{\nu+1/2}}.
\label{sdzd}%
\end{equation}
Next, we employ two integral relations related to the Hankel transform (see
\cite{Erdelyi_TIT_Vol2}, \cite{Oberhettinger_bessel}). Let $J_{\nu}$ denote
the Bessel function of first kind and order $\nu$ (see \cite[10.2.2]%
{NIST:DLMF}, \cite[p. 41]{Bessel_J_orig}). The first one reads (notation as in
(\ref{shorthand2})):%
\[
\frac{\operatorname*{e}\nolimits^{-s\left\Vert \mathbf{z}\right\Vert }%
}{\left(  \left\Vert \mathbf{z}\right\Vert +z_{d}\right)  ^{\nu}\left\Vert
\mathbf{z}\right\Vert }=\int_{0}^{\infty}\left(  \frac{x}{\omega}\right)
^{\nu+1/2}\frac{\operatorname*{e}\nolimits^{-\tilde{q}z_{d}}}{\tilde{q}\left(
\tilde{q}+s\right)  ^{\nu}}J_{\nu}\left(  x\omega\right)  \sqrt{x\omega}\ dx,
\]
see \cite[2.12.10.13]{PrudnikovVol2}. The second relation is taken from
\cite[p.31 (22)]{Erdelyi_TIT_Vol2}, \cite[(5.20)]{Oberhettinger_bessel},
\cite[2.12.10.10]{PrudnikovVol2}):%
\[
2\left(  \frac{s}{2\pi}\right)  ^{\nu+3/2}\frac{K_{\nu+1/2}\left(  sr\right)
}{r^{\nu+1/2}}=\frac{s}{\left(  2\pi\right)  ^{\nu+1}}\int_{0}^{\infty}\left(
\frac{x}{\omega}\right)  ^{\nu+1/2}\frac{\operatorname*{e}\nolimits^{-\tilde
{q}z_{d}}}{\tilde{q}}J_{\nu}\left(  x\omega\right)  \sqrt{x\omega}\ dx.
\]
We insert these relations into (\ref{sdzd}) and obtain after some
straightforward manipulations that (\ref{sdzd}) is equivalent to%
\begin{equation}
\underset{=:L}{\underbrace{\left(  s-\frac{\partial}{\partial z_{d}}\right)
^{\nu}\int_{0}^{\infty}\left(  \frac{x}{\omega}\right)  ^{\nu+1/2}%
\frac{\operatorname*{e}\nolimits^{-\tilde{q}z_{d}}}{\tilde{q}\left(  \tilde
{q}+s\right)  ^{\nu}}J_{\nu}\left(  x\omega\right)  \sqrt{x\omega}\ dx}}%
=\int_{0}^{\infty}\left(  \frac{x}{\omega}\right)  ^{\nu+1/2}\frac
{\operatorname*{e}\nolimits^{-\tilde{q}z_{d}}}{\tilde{q}}J_{\nu}\left(
x\omega\right)  \sqrt{x\omega}\ dx. \label{finalrelation}%
\end{equation}
We interchange the differentiation on the left-hand side with the integration
and make use of the simple dependence of the integrand on $z_{d}$ only through
the exponential factor, more precisely, we employ%
\[
\left(  s-\frac{\partial}{\partial z_{d}}\right)  ^{\nu}\operatorname*{e}%
\nolimits^{-\tilde{q}z_{d}}=\left(  s+\tilde{q}\right)  ^{\nu}%
\operatorname*{e}\nolimits^{-\tilde{q}z_{d}}.
\]
In this way, the left-hand side $L$ in (\ref{finalrelation}) equals
\[
L=\int_{0}^{\infty}\left(  \frac{x}{\omega}\right)  ^{\nu+1/2}\frac
{\operatorname*{e}\nolimits^{-\tilde{q}z_{d}}}{\tilde{q}}J_{\nu}\left(
x\omega\right)  \sqrt{x\omega}\ dx,
\]
and this is the right-hand side in (\ref{finalrelation}).%
\endproof

\textbf{Acknowledgement: }The first author gratefully acknowledges the
financial support by the Swiss National Science Foundation under grant no.
200020\_196995 and the second author the support by the Austrian Science Fund
(FWF) through the SFB F65.

We are grateful to our colleagues Prof.~Valery Smyshlaev and Prof.~Yury
Brychkov for interesting discussions concerning the explicit evaluation of
Sommerfeld-type integrals.\smallskip

\textbf{Conflict of interest statement: }This work does not have any conflicts
of interest.\smallskip

\textbf{Data Availability Statement:} Data sharing not applicable to this
article as no datasets were generated or analysed during the current study.

\bibliographystyle{abbrv}
\bibliography{nlailu}

\end{document}